\documentclass[a4paper,12pt]{article}
\usepackage{amsmath,amssymb,color}
\usepackage{latexsym}
\usepackage{graphicx}
\usepackage{float}
\newtheorem{theorem}{Theorem}

\begin{document}
\title{Fractional Factors, Component Factors and Isolated Vertex Conditions in Graphs}

\author{
Mikio Kano\footnote{Email: mikio.kano.math@vc.ibaraki.ac.jp;
Supported by JSPS KAKENHI Grant Number 16K05248.}\\
Ibaraki University, Hitachi, Ibaraki, Japan\\
Hongliang Lu\footnote{Email: luhongliang@mail.xjtu.edu.cn;
Supported by the National Natural
Science Foundation of China under grant No.11871391 and
Fundamental Research Funds for the Central Universities.}
\\ School of Mathematics and Statistics\\
Xi'an Jiaotong University,
Xi'an, Shaanxi, China\\
Qinglin Yu\footnote{Email: yu@tru.ca; Supported by the Discovery Grant of NSERC of Canada, and Shaanxi Hundred-Talents Program.}\\
School of Science\\
Xi'an Polytechnic University, Xi'an, Shaanxi, China\\
Department of Mathematics and Statistics\\
Thompson Rivers University, Kamloops, BC, Canada
}

\date{}

\maketitle
\begin{abstract}
For a graph $G = (V, E)$, a {\em fractional $[a, b]$-factor} is a real valued function $h:E(G)\to [0,1]$ that satisfies
$a \le ~ \sum_{e\in E_G(v)} h(e) ~ \le b$ for all $ v\in V(G)$, where $a$ and $b$ are real numbers and  $E_G(v)$ denotes the set of edges incident with $v$.
In this paper, we prove that the condition $\mathit{iso}(G-S) \le (k+\frac{1}{2})|S|$ is equivalent to the existence of fractional $[1,k+ \frac{1}{2}]$-factors, where ${\mathit{iso}}(G-S)$ denotes the number of isolated vertices in $G-S$. Using fractional factors as a tool, we construct component factors under the given isolated conditions. Namely, (i) a graph $G$ has a $\{P_2,C_3,P_5, \mathcal{T}(3)\}$-factor if and only if $\mathit{iso}(G-S) \le \frac{3}{2}|S|$ for all $S\subset V(G)$;
(ii) a graph $G$ has a $\{K_{1,1}, K_{1,2}, \ldots,$ $K_{1,k}, \mathcal{T}(2k+1)\}$-factor ($k\ge 2$)
if and only if $\mathit{iso}(G-S) \le (k+\frac{1}{2})|S|$ for all $S\subset V(G)$,
where $\mathcal{T}(3)$ and $\mathcal{T}(2k+1)$ are two special families of trees.

\end{abstract}

\bigskip
\noindent {\bf Keywords:} component factors, fractional factors, isolated vertices, isolated toughness, Tutte type condition.

%%%%%%%%%%%%%%%%%%%%%%%%%%%%%%%%%%%

\section{Introduction}

In this paper, we mainly consider finite simple graphs, which have neither loops nor multiple edges. A graph that has multiple edges but has no loops is referred as a {\em multigraph}. When defining notation and definitions,  we often referred a multigraph as a graph for convenience. Let $G$ be a graph with vertex set $V(G)$ and edge set $E(G)$. The number of vertices of $G$ is called its {\em order} and denoted by $|G|$. On the other hand, the number of edges in $G$ is called its {\em size} and denoted by $||G||$.

For a set $X$, the cardinality of $X$ is denoted by $|X|$ or $\#X$.
For a vertex $v$ of a graph $G$, the degree of $v$ in $G$ is denoted by $\deg_G(v)$. For two vertices $x$ and $y$ of $G$, an edge joining them is denoted by $xy$ or $yx$. We denote by $\mathit{Iso}(G)$ the set of isolated vertices of $G$, and  by ${\mathit{iso}}(G)$  the number of isolated vertices in $G$. Thus $\mathit{iso}(G-S)=|\mathit{Iso}(G)|$. For two disjoint vertex sets $X$ and $Y$ of $G$, the set of edges of $G$ joining $X$ to $Y$ is written as $E_G(X,Y)$ and $e_G(X,Y):=|E_G(X,Y)|$.
$P_n$ and $C_n$ are the path and the cycle of order $n$, respectively.
The set of non-negative integers is denoted by $\mathbf{Z^*}=\{0\} \cup \mathbf{Z^+}$.

Let $G$ be a graph, and $g,f:V(G) \to \mathbf{Z^*}$ be two integer-valued functions with $g\le f$, that is, $0\le g(x)\le f(x)$ for all $x\in V(G)$. Then a spanning subgraph $F$ of $G$ is called a $(g,f)$-{\it factor} of $G$ if $g(x) \le \deg_F(x) \le f(x)$ for all  $x\in V(G)$. For a vertex $v$ of $G$, let $E_G(v)$ denote the set of edges of $G$ incident with $v$.
For real-valued functions $g,f:V(G) \to \mathbf{R}$ with $g\le f$, a {\em fractional $(g,f)$-factor} is a function $h:E(G)\to [0,1]$ that satisfies the following condition:
\begin{align}
g(v) \le ~ \deg^h(v):=\sum_{e\in E_G(v)} h(e) ~ \le f(v) \quad \mbox{for all} \quad v\in V(G),
\end{align}
\noindent where $\deg^h(v)$ is called the {\em  $h$-degree of $v$} and $h(e)$ is a real number between 0 and 1 including 0 and 1. If the values of $h$ are 0 and 1 only, then a fractional $(g,f)$-factor becomes a $(g,f)$-factor.
Many results on fractional factors of graphs can be found in \cite{YuLiu}.

To study fractional factors, Yang, Ma and Liu \cite{YML} introduced a new parameter, {\it isolated toughness} of a graph $G$, denoted by $I(G)$, which is defined as
$$I(G)=\min \Big\{\frac{|S|}{iso(G-S)}\ :\ S\subseteq V(G),iso(G-S)\geq 2 \Big\} $$
if $G$ is not complete. Otherwise, $I(G)=\infty$. A graph $G$ is called \emph{isolated $t$-tough} if $I(G)\geq t$, where $t\geq 0$ is a real number.

For a set $\mathcal{S}$ of connected graphs, a spanning subgraph $F$ of $G$ is called an $\mathcal{S}${\em -factor} if each component of $F$ is isomorphic to an element of $\mathcal{S}$ (see (1) of Figure~1).
For a set $\mathbb{S}$ of positive integers, a spanning subgraph $F$ of $G$ is called an $\mathbb{S}${\em -factor} of $G$ if $\deg_F(x) \in \mathbb{S}$ for all vertices $x$ of $G$.
For an integer $k\ge 0$, the set of vertices of $G$ with degree $k$ is denoted by $V_k(G)$, namely,  $V_k(G)=\{v\in V(G): \deg_G(v)=k\}$.

For a tree $T$, the set of leaves is denoted by $\mathit{Leaf}(T)$, i.e., $V_1(T)$. An edge of $T$ incident with a leaf is called a {\em pendant edge}. In particular, the number of leaves of $T$ is equal to that of pendant edges of $T$.

\begin{figure}[H]
  \centering
  \includegraphics[width=10cm]{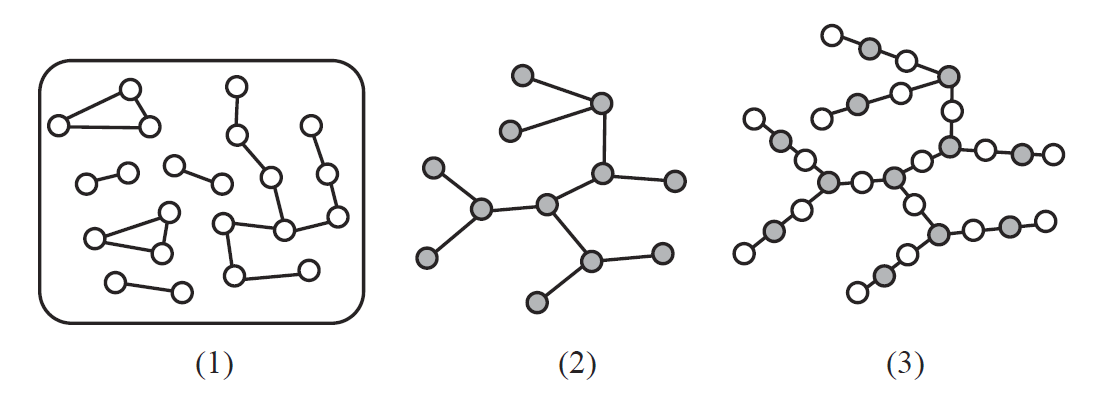}\\
  \vspace{-2mm}
  \caption{(1) A $\{P_2, C_3, P_5, \mathcal{T}(3)\}$-factor, which contains no $P_5$-components;  (2) A $\{1,3\}$-tree $R$; (3) The resulting $\{1,2,3\}$-tree $T_R$ obtained from $R$.}
 \label{fig-1}
\end{figure}

We define a special class of trees $\mathcal{T}(3)$: for any $\{1,3\}$-tree $R$ (i.e., every vertex has degree 1 or 3), a new tree $T_R$ is obtained from $R$ by inserting a new vertex of degree 2 into every edge of $R$, and by adding a new  pendant edge together its endpoint to every leaf of $R$ (see (2) and (3) of Figure~\ref{fig-1}). Then the tree $T_R$ is a $\{1,2,3\}$-tree having $||R||+|\mathit{Leaf}(R)|$ vertices of degree 2 and has the same number of leaves as $R$. Also there is one-to-one correspondence between $V_3(R)$ and $V_3(T_R)$. The collection of such $\{1,2,3\}$-trees $T_R$ generated from all $\{1,3\}$-trees $R$  is denoted by $\mathcal{T}(3)$. A more general class of trees, $\mathcal{T}(2k+1)$ ($k\ge 2$), will be defined in Section 2.
\medskip

Tutte \cite{Tutte1953} established a relationship between isolated 1-tough graphs and  $\{K_2,C_n: n \ge 3\}$-factors.

\begin{theorem}[Tutte \cite{Tutte1953}]\label{Tutte}
 A graph $G$ has a $\{K_2,C_n: n \ge 3\}$-factor if and only if
\[ iso(G-S) \le |S| \quad \mbox{for all} \quad S\subset V(G). \]
\end{theorem}

Amahashi and Kano \cite{AK1982} extended Theorem \ref{Tutte} and gave a characterization for isolated $1/k$-tough graphs in term of star factors.

\begin{theorem}[Amahashi and Kano \cite{AK1982}] Let $k\ge 2$ be an integer.
 A graph $G$ has a $\{K_{1,j} : 1\le j \le k \}$-factor if and only if
\[ iso(G-S) \le k|S| \quad \mbox{for all} \quad S\subset V(G). \]
\end{theorem}

Kano, Lu and Yu \cite{KLY2010} obtained a sufficient condition for isolated 2-tough graphs to have a components factor.
\begin{theorem}[Kano, Lu and Yu \cite{KLY2010}]  A graph $G$ has
a $\{K_{1,2}, K_{1,3},K_{5}\}$-factor if
\[ iso(G-S) \le \frac{|S|}{2} \quad \mbox{for all} \quad S\subset V(G). \]
\end{theorem}

Kano and Saito \cite{KS2012} as well as Zhang, Yan and Kano \cite{ZYK2015}  used isolated $k$-toughness to ensure the existence of special classes of component factors.

\begin{theorem}[Kano and Saito \cite{KS2012}] Let $k\ge 2$ be an integer.
 A graph $G$ has a $\{K_{1,j} : 1\le j \le 2k \}$-factor if
\[ iso(G-S) \le \frac{|S|}{k} \quad \mbox{for all} \quad S\subset V(G). \]
\end{theorem}

\begin{theorem}[Zhang, Yan and Kano \cite{ZYK2015}] Let $k\ge 2$ be an integer.
 A graph $G$ has a $\{K_{1,j}, K_{2k} : 1 \le j \le 2k-1\}$-factor if
\[ iso(G-S) \le \frac{|S|}{k} \quad \mbox{for all} \quad S\subset V(G). \]
\end{theorem}

In this paper, we carry on the investigations along the same direction mentioned above and obtain the factor characterizations of $I(G)= \frac{2}{3}$ and $I(G)= \frac{2}{2k+1}$ ($k \ge 2$):

\begin{theorem} \label{thm-1}
A graph $G$ has a $\{P_2,C_3,P_5, \mathcal{T}(3)\}$-factor if and only if
\begin{align}
\mathit{iso}(G-S) \le \frac{3}{2}|S| \quad \mbox{for all}\quad S\subset V(G).
\label{eq-17}
\end{align}
\end{theorem}

\begin{theorem} \label{thm-3}
Let $k\ge 2$ be an integer. Then a graph $G$ has a $\{K_{1,1}, K_{1,2}, \ldots,$ $K_{1,k}, \mathcal{T}(2k+1)\}$-factor if and only if
\begin{align}
\mathit{iso}(G-S) \le \Big( k+ \frac{1}{2}\Big) |S| \quad \mbox{for all}\quad S\subset V(G).
\label{eq-20}
\end{align}
\end{theorem}

\section{Proofs of Theorems \ref{thm-1} and \ref{thm-3}}

For a function $f:V(G) \to \mathbf{Z^*}$ and a vertex set $X$ of $G$, we write
\[ f(X) :=\sum_{x\in X}f(x), \quad \mbox{in particular}, \quad \deg_G(X)=\sum_{x\in X}\deg_G(x). \]

The tools for proving Theorems \ref{thm-1} and \ref{thm-3} are fractional factors. We first characterize the condition (\ref{eq-20}) in fractional $[1,k+\frac{1}{2}]$-factors (Theorem \ref{thm-15}), and then show that the minimal fractional factors are the desired component factors. In establishing the link between the condition (\ref{eq-20}) and fractional factors, we need the following theorem.

\begin{theorem}[Heinrich et al. \cite{HeinHellKirkLiu-90}, Anstee  \cite{An90}] Let $G$ be a multigraph and $g,f:V(G)\to \mathbf{Z^*}$ with $g< f$. Then $G$ has a $(g,f)$-factor if and only if
\begin{align}
g(T)-\deg_{G-S}(T) \le f(S) \quad \mbox{for all} \quad S\subset V(G),
\label{eq-10}
\end{align}
where $T=\{v\in V(G)-S:  \deg_{G-S}(v)<g(v)\}$.
\label{thm-10}
\end{theorem}

If two functions $g,f:V(G)\to \mathbf{R}$ take the constant real values $g(x)=a$ and $f(x)=b$ for every vertex $x$, then fractional $(g,f)$-factors are called {\em fractional $[a,b]$-factors} in short.

\begin{theorem} Let $k\ge 1$ be an integer and $G$ be a graph. Then $G$ has a fractional $[1,k+\frac{1}{2}]$-factor $h$ with values in $\{0,\frac{1}{2},1\}$ if and only if
 \begin{align}
\mathit{iso}(G-S) \le \Big(k+ \frac{1}{2}\Big) |S| \quad \mbox{for all}\quad S\subset V(G). \label{eq-30}
\end{align}
\label{thm-15}
\end{theorem}
\vspace{-1cm}

\noindent {\em Proof.} Assume that $G$ satisfies (\ref{eq-30}). Let $G^*$ denote the multigraph obtained from $G$ by replacing each edge $e$ of $G$ by two parallel edges $e(1)$ and $e(2)$. Then $V(G^*)=V(G)$, and $\deg_{G^*}(v)=2\deg_G(v)$ for every  $v\in V(G^*)$. Define two functions $g,f:V(G^*)\to \mathbf{Z^*}$ as
\[ g(x)=2 \quad \mbox{and} \quad f(x)=2k+1 \quad \mbox{for all}\quad x\in V(G^*).\]
Then $g<f$, and for every  $S\subset V(G^*)$, we have
\begin{align*}
 T & =\{v\in V(G^*)-S: \deg_{G^*-S}(v)<g(v)=2\} \\
& =\{v\in V(G^*)-S: \deg_{G^*-S}(v)=0\}.  \\
& = \mathit{Iso}(G-S).
\end{align*}
Thus it follows from the above equality and (\ref{eq-30}) that
\begin{align*}
g(T)-\deg_{G^*-S}(T)  & =2\cdot \mathit{iso}(G-S) -0  \\
& \le 2\cdot \Big(k+ \frac{1}{2}\Big) |S| =(2k+1)|S|=f(S).
\end{align*}
Hence by Theorem~\ref{thm-10}, $G^*$ has a  $(g,f)$-factor $F$. Now we construct a fractional $[1,k+\frac{1}{2}]$-factor $h:E(G)\to \{0,\frac{1}{2},1\}$ as follows: for an edge $e$ of $G$, (i) if $F$ contains both edges $e(1)$ and $e(2)$ of $G^*$, then define $h(e)=1$; (ii) if $F$ contains exactly one of $e(1)$ and $e(2)$, then define $h(e)=\frac{1}{2}$; (iii) otherwise, let $h(e)=0$. It is easy to see that $h$ is the desired fractional $[1,k+\frac{1}{2}]$-factor with values in $\{0,\frac{1}{2},1\}$.

\medskip

Next assume that $G$ has a fractional $[1,k+\frac{1}{2}]$-factor $h$ with values in $ \{0,\frac{1}{2},1\}$.
Let $S\subset V(G)$, and $F$ be the spanning subgraph of $G$ induced by $\{e\in E(G): h(e)=\frac{1}{2}~\mbox{or}~ 1\}$.  Clearly, the neighbors of each isolated vertex $u$ of $G-S$ are contained in $S$ and $\deg^h(u)\ge 1$, thus we have
\begin{align*}
\mathit{iso}(G-S) & \le  \sum_{e\in E_F(\mathit{Iso}(G-S),S)} h(e)  \\
& \le \sum_{x\in S} \deg^h(x) \le \Big( k+\frac{1}{2} \Big)|S|.
\end{align*}
Hence $\mathit{iso}(G-S) \le (k+\frac{1}{2})|S|$, i.e., (\ref{eq-30}) holds. \hfill $\Box$

\bigskip

\noindent
{\em Proof of Theorem \ref{thm-1}.} We first show that every tree $T \in \mathcal{T}(3)$ satisfies the condition (\ref{eq-17}). Define a function $h:E(T)\to \{\frac{1}{2},1\}$ as follows: for every pendant edge $e_1$ of $T$, let $h(e_1)=1$ and for any other edge $e_2$, let $h(e_2)=\frac{1}{2}$. Since $T$ is a $\{1,2,3\}$-tree and no pendant edge is incident with a vertex of degree 3, $h$ is a fractional $[1,\frac{3}{2}]$-factor. Hence, by Theorem~\ref{thm-15}, $T$ satisfies the condition (\ref{eq-17}).

Assume that $G$ has a $\{P_2,C_3,P_5,\mathcal{T}(3)\}$-factor $F$.
Let $D_1,D_2, \ldots, D_m$ be the components of $F$. Then each $D_i$ is $P_2$, $C_3$, $P_5$ or a tree of $\mathcal{T}(3)$, and thus $\mathit{iso}(D_i-X_i)\le \frac{3}{2}|X_i|$ for every $X_i\subset V(D_i)$.
Then for every  $S\subset V(G)$, we have
\begin{align*}
 \mathit{iso}(G-S) & \le \mathit{iso}(F-S)=\sum_ {i=1}^m \mathit{iso}(D_i-S\cap V(D_i)) \\
& \le \sum_{i=1}^m \frac{3}{2}|S\cap V(D_i)| = \frac{3}{2}|S|.
\end{align*}
Hence the necessity is proved.

\bigskip

Next we prove the sufficiency. By Theorem~\ref{thm-15}, $G$ has a fractional $[1,\frac{3}{2}]$-factor $h$ with values in  $\{0,\frac{1}{2},1\}$. We call an edge $e_1$ with $h(e_1)=1$ a {\em red edge} and an edge $e_2$ with $h(e_2)=\frac{1}{2}$ a {\em blue edge}. Let $F$ be the subgraph of $G$ induced by the set of all red and blue edges. Namely, $F$ is obtained from $G$ by removing all the edges $e_3$ with $h(e_3)=0$. Since $h$ is a fractional $[1,\frac{3}{2}]$-factor, $F$ is a spanning subgraph of $G$, and for every vertex $v$ of $G$, one of the following two statements holds:

\medskip
 (i) no red edge is incident with $v$ and two or three blue edges are incident with $v$;  or

 (ii) exactly one red edge is incident with $v$ and at most one blue edge is incident with $v$.

\medskip \noindent
Choose a fractional $[1,\frac{3}{2}]$-factor $h$ so that the number of edges in $F$ is as small as possible. For a convenience, we also call $F$ a fractional $[1,\frac{3}{2}]$-factor.
For every vertex $v$, we call the number of edges of $F$ incident with $v$ {\em the degree} of $v$ and denote it by $\deg_F(v)$. It is clear that $\deg_F(v)\in \{1,2,3\}$. On the other hand, the degree $\deg^h(v)$ of $v$ in $h$ is 1 or $\frac{3}{2}$.

\medskip \noindent
{\bf Claim~1.} {\em Every cycle of $F$ is a $C_3$-component of $F$.}
\medskip

\noindent
{\em Proof.} Assume that $F$ contains a cycle $C$. First assume that $C$ is of even order. Take a perfect matching $M$ of $C$, and recolor all the edges of $M$ red, and remove all the edges in $C-M$. Then the resulting subgraph is a new fractional $[1,\frac{3}{2}]$-factor with red and blue edges, but its size is smaller than $F$, which contradicts the choice of $F$. Hence $C$ is of odd order.

Assume that $C$ has two adjacent vertices $u_1$ and $u_2$ with degree 3 in $F$. Then $F-u_1u_2$ is a new fractional $[1,\frac{3}{2}]$-factor with fewer edges than $F$, which contradicts to the choice of $F$. Hence if a vertex $v$ of $C$  has degree 3 in $F$, then the two neighbors of $v$ in $C$ have degree 2 in $F$. Assume that $C$ has a vertex $v$ with $\deg_F(v)= 3$. Let $u_1$ and $u_2$ be the two neighbors of $v$ in $C$. Take a perfect matching $M$ of $C-v$. Recolor the edges of $M$ red, and remove all the edges of $(C-v)-M$ and  $vu_2$. Since $vu_1$ and $vu_2$ are both blue edges, we obtain a new fractional $[1,\frac{3}{2}]$-factor with fewer edges than $F$, which is a contradiction. Hence $C$ is a component of $F$.

Moreover, it is easy to see that an odd cycle of $C$ order at least 5 has a $\{P_2,P_5\}$-factor $F_C$. Remove all the edges of $C$ not contained in $F_C$, recolor the edges contained in $P_2$ of $F_C$ red, and two pendant edges of $P_5$ of $F_C$ red and the remaining two edges of $P_5$ of $F_C$ blue. Then we obtain a new fractional $[1,\frac{3}{2}]$-factor with fewer edges than $F$, a contradiction. Therefore every cycle contained in $F$ is $C_3$. Consequently Claim 1 is proved.

\medskip

For the simplicity of statements, from now on, we will use ``another fractional $[1,\frac{3}{2}]$-factor" to replace the phase ``a new fractional $[1,\frac{3}{2}]$-factor with fewer edges than $F$".

\medskip \noindent
{\bf Claim~2.} {\em Every non-cycle component of $F$ is $P_2$, $P_5$ or a tree of $\mathcal{T}(3)$.}
\medskip

\noindent
{\em Proof.} Let $x$ and $y$ be two vertices of degree 3 in $F$ such that they are adjacent or connected by a path whose all inner vertices have degree 2 in $F$. If $x$ and $y$ are adjacent in $F$, then $F-xy$ is another fractional $[1,\frac{3}{2}]$-factor, a contradiction. Assume that $x$ and $y$ are connected by a path $(x,u_1,u_2,\ldots, u_n,y)$ of length at least 3 (i.e., $n\ge 2$) such that every $u_i$ has degree 2 in $F$. Then remove  $u_ny$, recolor $u_{n-1}u_n$ red, and recolor all remaining edges of the path blue. Then resulting  subgraph is another fractional $[1,\frac{3}{2}]$-factor, a contradiction.
Therefore,

\medskip

 (a) {\em if two vertices of $F$ with degree 3 are connected by a path in $F$ whose inner vertices have degree 2 in $F$, then the length of the path is 2.}

\medskip
Let $z$ be a leaf of $F$ and $x$ be a vertex of degree 3 in $F$. If $z$ and $x$ are adjacent,  then the edge $xz$ is red and so $\deg^h(x)\ge 2$, which is impossible. Hence $z$ and $x$ are not adjacent.  Assume that $z$ and $x$ are connected by a path $(z,u_1,u_2,\ldots, u_n,x)$ with $\deg_F(u_i)=2$ for every $i$. Then $zu_1$ is red.
First assume that  $n\ge 3$. Then remove $u_{n}x$, recolor $u_{n-1}u_n$ red, and recolor all the remaining edges of the path except $zu_1$ blue. Then the resulting subgraph is another fractional $[1,\frac{3}{2}]$-factor, a contradiction.  Next assume $n=1$. Then by removing $u_1x$, we obtain another fractional $[1,\frac{3}{2}]$-factor, a contradiction again.
Therefore,

\medskip
 (b) {\em if a leaf $z$ is contained in a component $D$ of $F$ with $\Delta(D)\ge 3$, then  there is a path $(z,u_1,u_2,x)$  in $F$ such that $\deg_F(u_i)=2$ and $\deg_F(x)= 3$.}

\medskip \noindent
Consequently, if a component $D$ of $F$ contains at least two vertices of degree 3, then by (a) and (b), $D$ is a tree of $\mathcal{T}(3)$. If $D$ has exactly one vertex of degree 3, then by the property (b), $D$ is also a tree of $\mathcal{T}(3)$.
If $D$ has no vertex of degree 3, then $D$ is a path. It is obvious that $P_3$ has no fractional $[1,\frac{3}{2}]$-factor, and so $D$ is not $P_3$. If $D$ is a path of even order, then $D$ has a $P_2$-factor, and so it contradicts to the minimality of $F$. If $D$ is a path $P_n$ of odd order with $n\ge 7$, then $D$ has a $\{P_2,P_5\}$-factor and thus contradicts to the minimality of $F$. Hence Claim~2 holds. %|~~$\Box$

\medskip
Clearly, Claims 1 and 2 imply the sufficiency. \hfill  $\Box$

%%\bigskip
%\begin{figure}
%  \centering
%  % Requires \usepackage{graphicx}
%  \includegraphics[width=10cm]{isolated-factor-2.png}\\
%  \caption{(1) A tree $R$ that satisfies (\ref{eq-15}) with $k=4$; (2) The tree $R-Leaf(R)$; (3) The tree $T_R$ obtained from $R$, where all red edges $e_1$ have  value $h(e_1)=1$ and all other edges $e_2$ have value $h(e_2)=\frac{1}{2}$.}
% \label{fig-2}
%\end{figure}

%\begin{figure}[htbp]
% \begin{center}
%% \includegraphics*{isolated-factor-2.pdf}
% \includegraphics*{isolated-factor-2.eps}
% \end{center}
% \caption{(1) A tree $R$ that satisfies (\ref{eq-15}) with $k=4$; (2) The tree $R-Leaf(R)$; (3) The tree $T_R$ obtained from $R$, where all red edges $e_1$ have  value $h(e_1)=1$ and all other edges $e_2$ have value $h(e_2)=\frac{1}{2}$.}
% \label{fig-2}
% \end{figure}

To state Theorem \ref{thm-3}, we need a new class of trees $\mathcal{T}(2k+1)$. Let $k\ge 2$  be an integer and $R$ be a tree that satisfies the
following conditions:  for every vertex $v \in V(R)$,
\begin{align}
\mbox{(i)} \quad & \deg_{R-Leaf(R)}(v) \in \{1,3,\ldots, 2k+1\}, \quad \mbox{and} \nonumber \\
\mbox{(ii)} \quad & 2 \cdot (\# \mbox{ of leaves adjacent to}~v) + \deg_{R-Leaf(R)}(v)~ \le 2k+1.
\label{eq-15}
\end{align}
\vspace{-2mm}
\begin{figure}[H]
  \centering
  \includegraphics[width=10cm]{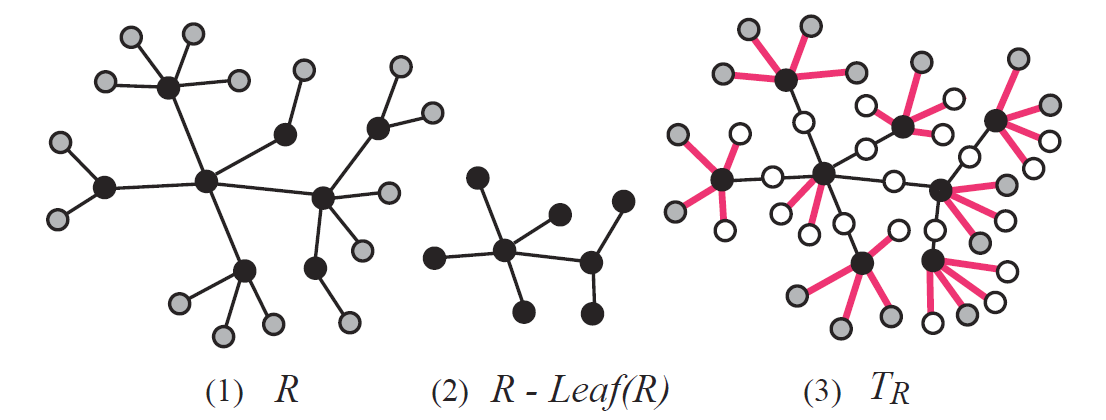}\\
  \caption{(1) A tree $R$ that satisfies (\ref{eq-15}) with $k=4$; (2) The tree $R-Leaf(R)$; (3) The tree $T_R$ obtained from $R$, where all red edges $e_1$ have  value $h(e_1)=1$ and all other edges $e_2$ have value $h(e_2)=\frac{1}{2}$.}
 \label{fig-2}
\end{figure}
\vspace{-2mm}
 (see  (1) and (2) of Figure~\ref{fig-2}). For such a tree $R$, we obtain a new tree $T_R$ as follows:

(iii)  insert a new vertex of degree 2 into each edge of $R-Leaf(R)$, and

(iv)  for each vertex $v$ of $R-Leaf(R)$ with $\deg_{R-Leaf(R)}(v)=2r+1<2k+1$, add  $ k-r-(\# \mbox{ of leaves adjacent to}~v)$ pendant edges together with their endpoints to $v$  (see (3) of Figure~\ref{fig-2} ).

\medskip  \noindent
Then the resulting tree $T_R$ has a fractional $[1,k+\frac{1}{2}]$-factor $h$ such that every pendant edge $e_1$ of $T_R$ has $h(e_1)=1$, all the other edges $e_2$ have $h(e_2)=\frac{1}{2}$. Moreover, this fractional  $[1,k+\frac{1}{2}]$-factor $h$ satisfies that every vertex $v$ of $R-Leaf(R)$  has $\deg^h(v)=k+\frac{1}{2}$, and every leaf $x$ of $T_R$ and every inserted vertex $y$ of degree 2 have $\deg^h(x)=\deg^h(y)=1$ (see  (3) of Figure~\ref{fig-2}).
The set of such trees $T_R$ for all trees $R$ satisfying  (\ref{eq-15})  is denoted by $\mathcal{T}(2k+1)$. Note that the construction of $\mathcal{T}(3)$ and that of $\mathcal{T}(2k+1)$ with $k\ge 2$ are similar, but adding pendant edges to some vertices of $R-Leaf(R)$ is not defined in the construction of $\mathcal{T}(3)$.

\bigskip

\noindent
{\em Proof of Theorem \ref{thm-3}.}  For any tree $T \in \mathcal{T}(2k+1)$, since $T$ has a fractional $[1,k+\frac{1}{2}]$-factor $h$ of values $\{\frac{1}{2}, 1\}$,
 $T$ satisfies (\ref{eq-20}) by Theorem~\ref{thm-15}.

Assume that $G$ has a $\{K_{1,1}, K_{1,2}, \ldots,$ $K_{1,k}, \mathcal{T}(2k+1)\}$-factor $F$.
Let $D_1,D_2, \ldots, D_m$ be the components of $F$. Then each $D_i$ is $K_{1,s}$ for some $1\le s \le k$, or a tree in $\mathcal{T}(2k+1)$. Thus $\mathit{iso}(D_i-X_i)\le (k+\frac{1}{2})|X_i|$ for every $X_i\subset V(D_i)$. Then, for any  $S\subset V(G)$, we have
\begin{align*}
 \mathit{iso}(G-S) & \le \mathit{iso}(F-S)=\sum_ {i=1}^m \mathit{iso}(D_i-S\cap V(D_i)) \\
& \le \sum_{i=1}^m \Big( k+ \frac{1}{2}\Big) |S\cap V(D_i)| = \Big( k+ \frac{1}{2}\Big) |S|.
\end{align*}
Hence the necessity is proved.

\medskip
Next we prove the sufficiency. Assume that $G$ satisfies (\ref{eq-20}).
By Theorem~\ref{thm-15}, $G$ has a fractional $[1,k+\frac{1}{2}]$-factor $h$ with values in  $\{0,\frac{1}{2},1\}$. We call an edge $e_1$ with $h(e_1)=1$ a {\em red edge} and an edge $e_2$ with $h(e_2)=\frac{1}{2}$ a {\em blue edge}. Let $F$ be the subgraph of $G$ induced by the set of all red and blue edges. Namely, $F$ is obtained from $G$ by removing all the edges $e_3$ with $h(e_3)=0$, and since $h$ is a fractional $[1,k+\frac{1}{2}]$-factor, $F$ is a spanning subgraph.
Choose a fractional $[1,k+\frac{1}{2}]$-factor $h$ of $G$ so that the number of edges in $F$ is as small as possible.

For a vertex $v$, we call the number of edges of $F$ incident with $v$ {\em the degree} of $v$ in $F$ and denote it by $\deg_F(v)$. It is clear that $1\le \deg_F(v) \le 2k+1$. On the other hand, the $h$-degree of $v$ can be expressed as $\deg^h(v)=1+\frac{1}{2}t$ for some integer $t$, $0\le t\le 2k-1$.

\medskip \noindent
{\bf Claim~1.} {\em $F$ contains no cycle, i.e., $F$ is a forest.}
\medskip

\noindent
{\em Proof.} Suppose that $F$ contains a cycle $C$. First assume that $C$ is of even order. Take a perfect matching $M$ of $C$, and recolor all edges of $M$ red, and remove all edges in $C-M$. Then the resulting subgraph is a new fractional $[1,k+\frac{1}{2}]$-factor of $G$, but its size is smaller than $F$, a contradiction. Hence $C$ is of odd order.

As in the proof of Theorem \ref{thm-1}, we will replace ``a new fractional $[1,k+\frac{1}{2}]$-factor with less edges than $F$" by ``another fractional factor" for the simplicity in the rest of the proof.

Assume that $C$ has two adjacent vertices $v_1$ and $v_2$ with degree at least 3 in $F$. Then $F-v_1v_2$ is another fractional factor, a contradiction. Hence if a vertex $v$ of $C$  has degree at least 3 in $F$, then the two neighbors of $v$ in $C$ have degree 2 in $F$. Assume that $C$ has a vertex $v$ with $\deg_F(v)\ge 3$. Let $u_1$ and $u_2$ be the two neighbors of $v$ in $C$. Take a perfect matching $M$ of $C-v$, and recolor all edges of $M$ red, and remove all edges of $(C-v)-M$ and the edge $vu_2$. Then we obtain another fractional factor,  a contradiction. Hence $C$ is a component of $F$.

It is easy to see that $C$ has a $\{P_2=K_{1,1}, P_3=K_{1,2}\}$-factor $F_C$.  We recolor all edges of $P_2$-components and $P_3$-components of $F_C$ red, and remove all other edges of $C$. Then we obtain another fractional factor, a contradiction. Therefore $F$ has no cycles,
and the claim is proved.

\medskip \noindent
{\bf Claim~2.} {\em Let $x$ and $y$ be two vertices of degree at least 3 in $F$. Then $x$ and $y$ are not adjacent in $F$. If $x$ and $y$ are connected by a path whose all inner vertices have degree 2 in $F$, then the length of the path is 2, and $\deg^h(x)=\deg^h(y)=k+\frac{1}{2}$ and the two edges in the path are blue edges.}

\medskip \noindent
{\em Proof.}  If $x$ and $y$ are adjacent in $F$, then $F-xy$ is another fractional factor,  a contradiction.
Assume that $x$ and $y$ are connected by a path $(x,u_1,u_2,\ldots, u_n,y)$ with $\deg_F(u_i)=2$ ($1\le i \le n$). If $n\ge 2$, then by removing  $u_ny$ and recoloring $u_{n-1}u_n$ red and all the remaining edges of the path blue, the resulting subgraph is another fractional factor, a contradiction. Therefore $n=1$, and the path is $(x,u_1,y)$.

If $\deg^h(x)<k+\frac{1}{2}$, then by removing $u_1y$ and recoloring $xu_1$ red, we obtain another fractional factor,  a contradiction. Hence  $\deg^h(x)=\deg^h(y)=k+\frac{1}{2}$ by the symmetry. If $xu_1$ is red, then removing $u_1y$ we obtain another fractional factor, a contradiction. Therefore $xu_1$ and $yu_1$ are blue edges, and the claim holds.

\medskip \noindent
{\bf Claim~3.} {\em A leaf $z$ in $F$ is either contained in a star component, or
adjacent to a vertex $x$ with $\deg_F(x)\ge 3$ and $\deg^h(x)=k+\frac{1}{2}$.}

\medskip \noindent
{\em Proof.}  Let $z$ be a leaf of $F$, and $D$ be the component of $F$ containing $z$.
Assume that $D$ is a path $(z,u_1,u_2, \ldots,$  $u_n,y)$ such that $\deg_F(u_i) = 2$ ($1\le i \le n$) and $\deg_F(y)=1$. If $n\ge 2$, then by removing $u_1u_2$ and recoloring $u_2u_3$ red, we obtain another fractional factor, a contradiction. If $n=1$, then $D=P_3=K_{1,2}$, which is a star.

Next assume that $F$ contains a path $(z,u_1,u_2, \ldots,$  $u_n,x)$ such that $\deg_F(u_i) = 2$ ($1\le i \le n$) and $x$ has degree at least 3 in $F$. If $n\ge 2$, then by removing $u_nx$, and recoloring $u_{n-1}u_n$ and $zu_1$ red, and recoloring all other remaining edges (if any) of the path blue, we obtain
another fractional factor, a contradiction. If $n=1$, then $F-u_1x$ is another fractional $[1,k+\frac{1}{2}]$-factor, a contradiction. Therefore $z$ and $x$ are adjacent.

Moreover, if $D$ contains exactly one vertex $x$ of degree at least 3, by the same argument given above, we see that every leaf of $D$ is adjacent to $x$ and thus $D$ is a star; otherwise, $D$ contains another vertex $y$ of degree at least 3, then by Claim~2, $\deg^h(x)=k+\frac{1}{2}$.
Thus the claim is proved.

\medskip \noindent
{\bf Claim~4.} (i) {\em If $uv$ is an edge such that $\deg_F(u)\ge 3$ and $\deg_F(v)=2$, then  $\deg^h(u)=k+\frac{1}{2}$ and $uv$ is a blue edge.} (ii) {\em If $xy$ is a red edge, then one of $x$ and $y$ is a leaf of $F$.}
\medskip

 \noindent
{\em Proof.} Let  $uv$ be an edge of $F$ such that $\deg_F(u)\ge 3$, $\deg_F(v)=2$ and $\deg^h(u)<k+\frac{1}{2}$. Let $z_1$ be a vertex adjacent to $v$. If $\deg_F(z_1) \ge 3$, then by Claim~2, we have $\deg^h(u)=k+\frac{1}{2}$, a contradiction. If $\deg_F(z_1)=1$, then it contradicts to Claim~3.
Hence $\deg_F(z_1)= 2$. By removing $uv$ and recoloring $vz_1$ red, we obtain another fractional factor, a contradiction. Hence if $uv$ is an edge with $\deg_F(u)\ge 3$ and $\deg_F(v)=2$, then $\deg^h(u)=k+\frac{1}{2}$. If the statement (ii) is true, then $uv$ is a blue edge, and so in order to show that $uv$ is a blue edge, it suffices to prove (ii).

Assume to the contrary that $xy$ is a red edge and none of $x$ and $y$ is a leaf. If $\deg_F(x)\ge 3$ and $\deg_F(y)\ge 3$, then $F-xy$ is another fractional factor, a contradiction. Hence we may assume that $\deg_F(x)= 2$ and $\deg_F(y)\ge 2$ by the symmetry. Let $z_1$ be another vertex adjacent to $x$. If $\deg_F(z_1)\ge 3$,
then $F-xz_1$ is another fractional factor since $xy$ is a red edge, a contradiction. By Claim~3, $\deg_F(z_1)=2$. From Claim~2, $F$ has no path connecting two vertices of degree at least 3 and passing through $(z_1,x,y)$. Hence we may assume that there is a path $(z,u_1, \ldots, u_n=z_1,x,y)$ with $\deg_F( z)=1$, $\deg_F(u_i)=2$ ($1\le i \le n$) and $\deg_F(x)=2$.
By Claim~3, $u_1=x$, which is a contradiction since $u_n=z_1\ne x$. Therefore Claim~4 holds.

\medskip \noindent
{\bf Claim~5.} {\em Every component $D$ of $F$ is isomorphic to one in $\{K_{1,1}, K_{1,2}, \ldots,$ $K_{1,k}, \mathcal{T}(2k+1)\}$.}
\medskip

\medskip \noindent
{\em Proof.} If $D$ is a star with center $u$, then $\deg^h(u)\le k+\frac{1}{2}$, which implies $\deg_D(u)\le k$ since every pendant edge of $D$ is red. So we may assume that $D$ is not a star, which implies that $D$ has at least two vertices of degree at least 3 by Claim~3.

Let $v$ be a vertex of degree at least 3. If $xv$ is a red edge, then $x$ is a leaf and vice versa. If $vy$ is a blue edge, then $\deg_F(y)=2$ and vice versa. Consequently, $D$ is a tree in $\mathcal{T}(2k+1)$.

\medskip

With Claim~5, we complete the proof. \hfill $\Box$

\section*{Acknowledgement}
The research was carried out during the first and the third authors' visits to Xi'an Jiaotong University. The support and hospitality from the university were greatly appreciated.
%We are grateful to an anonymous referee who provided with us many detailed
%suggestions on representation and on English.

\end{document}